\documentclass[12pt]{article}
\usepackage{amssymb,amsmath,amsfonts,amsthm}
\newcommand{\D}{{\Delta}}

\newtheorem{Th}{Theorem}
\newtheorem{conj}{Conjecture}
\newtheorem{lem}[Th]{Lemma}

\newtheorem{cor}[Th]{Corollary}

\newtheorem{obs}[Th]{Observation}
\title{Equitable Colorings of Planar Graphs without Short Cycles}

\author {Keaitsuda Nakprasit\footnote{Corresponding Author} \\ {\small\em Department of Mathematics, Faculty of Science, Khon Kaen University, 40002, Thailand }\\  
{\small\em E-mail address: kmaneeruk@hotmail.com} 
\and Kittikorn Nakprasit \\ 
{\small\em Department of Mathematics, Faculty of Science, Khon Kaen University, 40002, Thailand }\\
{\small\em E-mail address: kitnak@hotmail.com}}
\date{}

\begin {document}
\maketitle

\begin{abstract}
An \emph{equitable coloring} of a graph is a proper vertex
coloring such that the sizes of every two color classes differ by
at most $1.$ Chen, Lih, and Wu conjectured that
every connected graph $G$ with maximum degree $\Delta \geq 2$
has an equitable coloring with $\Delta$ colors, except when $G$ is a complete graph
or an odd cycle or $\Delta$ is odd and $G=K_{\Delta,\Delta}.$
Nakprasit proved the conjecture holds for planar graphs
with maximum degree at least $9.$
Zhu and Bu proved that the conjecture holds
for every $C_3$-free planar graph with maximum degree at least $8$ and
for every planar graph without $C_4$ and $C_5$ with maximum degree at least $7.$

In this paper, we prove that the conjecture holds
for planar graphs in various settings, especially
for every $C_3$-free planar graph with maximum degree at least $6$ and
for every planar graph without $C_4$ with maximum degree at least $7,$
which improve or generalize results on equitable coloring by Zhu and Bu.
Moreover, we prove that the conjecture holds
for every planar graph of girth at least $6$ with maximum degree at least $5.$

{\bf Key Words:} Equitable coloring; Planar graph; Cycle; Girth
\end{abstract}

\section{Introduction}

Throughout this paper, all graphs are finite, undirected, and simple.
We use  $V(G),$ $|G|,$ $E(G),$ $e(G),$ $\Delta (G),$ and $\delta(G),$ respectively, to
denote vertex set, order, edge set, size, maximum degree,
and minimum degree of a graph $G.$
We write $xy \in E(G)$ if $x$ and $y$ are adjacent.
The graph obtained by deleting an edge $xy$ from $G$ is denoted by $G - \{xy\}.$
For any vertex $v$ in $V(G),$ let $N_G(v)$ be the set of
all neighbors of $v$ in $G.$ The \emph{degree of $v$}, denoted by
$d_G(v),$ is equal to $|N_G(v)|.$
We use $d(v)$ instead of $d_G(v)$ if no confusion arises.
For disjoint subsets $U$ and $W$ of $V(G),$
the number of edges with one end in $U$ and another in $W$
is denoted by $e(U,W).$
We use $G[U]$ to denote the subgraph of $G$ induced by $U.$

An \emph{equitable $k$-coloring} of a graph is a proper vertex $k$-coloring 
such that the sizes of every two color classes differ by at most $1.$ We say that $G$ is
\emph {equitably $k$-colorable} if $G$ has an equitable $k$-coloring.

It is known \cite{GareyJohnson} that determining if a planar graph with
maximum degree $4$ is $3$-colorable is NP-complete. For a given
$n$-vertex planar graph $G$ with maximum degree $4,$ let $G'$
be a graph obtained from $G$ by adding $2n$ isolated vertices. Then $G$ is
$3$-colorable if and only if $G'$ is equitably $3$-colorable.
Thus, finding the minimum number of colors need to color a graph equitably
even for a planar graph is an NP-complete problem.

Hajnal and Szemer\'edi~\cite{HS} settled a conjecture of Erd\H os by
proving that  every graph $G$ with maximum degree at most $\Delta$
has an equitable $k$-coloring for every $k\geq 1+\Delta.$ In its
`complementary' form this result concerns decompositions of a
sufficiently dense graph into cliques of equal size. This result is
now known as Hajnal and Szemer\'edi Theorem. Later, Kierstead and
Kostochka~\cite{KK08} gave a simpler proof of Hajnal and Szemer\'edi
Theorem in the direct form of equitable coloring. The bound of the
Hajnal-Sz{e}mer\' edi theorem is sharp, but it can be improved for
some important classes of graphs. In fact, Chen, Lih, and
Wu~\cite{CLW94} put forth the following conjecture.

\begin{conj}
Every connected graph $G$ with
maximum degree $\Delta\geq 2$ has an equitable coloring with
$\Delta$ colors, except when $G$ is a complete graph or an odd
cycle or $\Delta$ is odd and $G=K_{\Delta,\Delta}.$
\end{conj}

Lih and Wu~\cite{LW} proved the conjecture for bipartite graphs.
Meyer \cite{M} proved that every forest with maximum degree $\Delta$
has an equitable $k$-coloring for each
$k \geq 1+\lceil \Delta/2\rceil $ colors.
This result implies conjecture holds for forests.
The bound of Meyer is attained at the complete bipartite $K_{1,m}:$ 
in every proper coloring of $K_{1,m},$ the center vertex forms a color class,
and hence the remaining vertices need at least $m/2$ colors.
Yap and Zhang~\cite{YZ1} proved that the conjecture holds for
outerplanar graphs.
Later Kostochka~\cite{Ko} extended the result for outerplanar graphs
by proving that every outerplanar graph
with maximum degree $\Delta$ has an equitable $k$-coloring for each
$k \geq 1+\lceil \Delta/2\rceil.$
Again this bound is sharp.

In~\cite{ZY98},
Zhang and Yap essentially proved the conjecture holds for planar graphs
with maximum degree at least $13.$
Later Nakprasit~\cite{Nak12} extended the result to all planar graphs
with maximum degree at least  $9.$

Other studies focused on planar graphs without some restricted cycles.
Li and Bu~\cite{LiBu09} proved that the conjecture holds for every planar graph without $C_4$ and $C_6$
with maximum degree at least $6.$
Zhu and Bu~\cite{ZhuBu08} proved that the conjecture holds for every $C_3$-free planar graph
with maximum degree at least $8$ and for every planar graph without $C_4$ and $C_5$ with
maximum degree at least $7.$
Tan \cite{Tan10} proved that the conjecture holds for every planar graph without $C_4$
with maximum degree at least  $7.$ Unfortunately the proof contains some flaws.

In this paper, we prove that each graph $G$ in various settings has 
an equitably $m$-colorable such that $m \leq \D.$  
Especially we  prove that the conjecture holds
for planar graphs in various settings, especially 
for every $C_3$-free planar graph with maximum degree at least $6$ and
for every planar graph without $C_4$ with maximum degree at least $7,$
which improve or generalize  results on equitable coloring by Zhu and Bu \cite{ZhuBu08}.
Moreover, we prove that the conjecture holds 
for every planar graph of girth at least $6$ with maximum degree at least $5.$

\section{Preliminaries}
Many proofs in this paper involve edge-minimal planar graph that
is not equitably $m$-colorable. The minimality is on inclusion, that is, any spanning subgraph with fewer edges is
equitably $m$-colorable. In this section, we describe
some properties of such graph that appear
recurrently in later arguments. The following fact about planar
graphs in general is well-known and can be found in standard texts about graph
theory such as~\cite{west}.

\begin{lem}\label{edge}
Every planar graph $G$ of order $n$ and girth $g$ has $e(G) \leq (g/(g-2))(n-2).$
Especially, a $C_3$-free planar graph $G$ has $e(G) \leq 2n-4$ and $\delta(G) \leq 3.$
\end{lem}

Let $G$ be an edge-minimal $C_3$-free planar graph that is not equitably $m$-colorable
with $|G|=mt,$ where $t$ is an integer. As $G$ is
planar and without $C_3,$ a graph $G$ has an edge $xy$ where $d(x) = \delta \leq 3.$ By
edge-minimality of $G,$ the graph $G - \{xy\}$ has an
equitable $m$-coloring $\phi$ having color classes $V'_1, V'_2, \ldots, V'_m.$ 
It suffices to consider only the case that $x, y \in V'_1.$
Choose $x \in V'_1$ such that $x$ has degree $\delta$ and order $V'_1, V'_2, \ldots ,V'_{\delta}$
in a way that $N(x) \subseteq V'_1 \cup V'_2 \cup \cdots \cup V'_{\delta}.$
Define  $V_1 = V'_1 - \{x\}$ and $V_i = V'_i$ for  $1 \leq i \leq m.$

We define $\mathcal R$ recursively. Let $V_1 \in \mathcal R$ and
$V_j \in \mathcal R$ if there exists a vertex in $V_j$ which has no
neighbors in $V_i$ for some $V_i \in \mathcal R.$ Let $r=|\mathcal
R|.$ Let $A$ and $B$ denote $\bigcup_{V_i \in \mathcal R} V_i$ and
$V(G) - A,$ respectively. Furthermore, we let $A'$ denote $A \cup \{x\}$ 
and $B'$ denote $B - \{x\}.$ From definitions of
$\mathcal R$ and $B, e(V_i,\{u\}) \geq 1$ for each $V_i \in \mathcal
R$ and $u \in B.$ Consequently $e(A,B)\geq r[(m-r)t+1]$ and
$e(A',B') \geq r(m-r)t.$

Suppose that there is $k$ such that $k \geq \delta +1$ and $V_k \in \mathcal R.$  
By definition of $\mathcal R,$ there exist $u_1 \in V_{i_1},u_2 \in
V_{i_2},\ldots,  u_s \in V_{i_s}, u_{i_{s+1}} \in V_{i_{s+1}}= V_k$
such that $e(V_1,\{u_1\}) = e(V_{i_1},\{u_2\}) = \cdots =
e(V_{i_s},\{u_{s+1}\}) = 0.$ Letting $W_1 = V_1 \cup \{u_1\},
W_{i_1} = (V_{i_1} \cup \{u_2\}) -\{u_1\}, \ldots, W_{i_s} =
(V_{i_s} \cup \{u_{s+1}\}) - \{u_s\},$ and  $W_k = (V_k\cup
\{x\}) - \{u_{s+1}\},$ otherwise $W_i=V_i,$ we get an
equitable $m$-coloring of $G.$  This contradicts to the fact that
$G$ is a counterexample.

Thus, in case of $C_3$-free planar graph,
we assume $\mathcal R \subseteq \{V_1, V_2, \ldots, V_\delta\}$
where $\delta \leq 3$ is the minimum degree of non-isolated vertices.

We summarize our observations here.

\begin{obs}\label{obs01} If $G$ is an edge-minimal $C_3$-free planar graph that
is not equitably $m$-colorable with order $mt,$ where $t$ is an integer, then we may assume\\
 \noindent (i) $\mathcal R \subseteq \{V_1, V_2,\ldots, V_{\delta}\}$ where $\delta \leq 3$
is the minimum degree of non-isolated vertices; \\
 \noindent (ii)  $e(u, V_i) \geq 1$ for each $u \in B$ and $V_i \in \mathcal R;$\\
 \noindent (iii) $e(A,B) \geq r[(m-r)t+1]$ and $e(A',B') \geq r(m-r)t.$
 \end{obs}

\section{Results on $C_3$-free Planar Graphs}

First, we introduce some useful tools and notation that will be used later. 

\begin{Th}\label{Grot} ~\cite{Grotzsch} (Gr\"{o}tzsch, 1959)
If $G$ is a $C_3$-free planar graph, then $G$ is $3$-colorable. 
\end{Th}

\begin{lem}\label{lem02}
Let  $m$ be a fixed integer with $m \geq 1.$ Suppose that any $C_3$-free
planar graph of order $mt$ with maximum degree at most $\D$ is
equitably $m$-colorable for any integer $t \geq k.$ Then any
$C_3$-free planar graph with order at least $kt$ and maximum degree
at most $\D$ is also equitably $m$-colorable.
\end{lem}

\noindent \textbf{Proof.} Suppose that any $C_3$-free planar graph
of order $mt$ with maximum degree at most $\D$ is equitably
$m$-colorable for any integer $t \geq k.$ Consider a $C_3$-free
planar graph $G$ of order $mt + r$ where $1\leq r \leq m-1$ and
$t\geq k.$ If $r=m-1$ or $m-2,$ then $G \cup K_{m-r}$ is equitably
$m$-colorable by hypothesis. Thus also is $G.$ Consider $r\leq m-3.$
Let $x$ be a vertex with minimum degree $d.$ We assume that 
$G - \{x\}$ is equitably $m$-colorable to use induction on $r.$
Thus the coloring of $G - \{x\}$ has $r+1$ color classes
with size $t-1.$ Since there are at most $d$ forbidden colors for
$x$ where $d\leq 3,$ we can add $x$ to a color class of size $t-1$ to
form an equitable $m$-coloring of $G.$ This completes the proof \qed

\begin{lem}\label{lem02a}~\cite{CLW94}
	If $G$ is a graph with maximum degree $\D \geq |G|/2,$ then $G$ is equitably $\D$-colorable.
\end{lem}

\begin{obs}\label{obs02}  By Lemmas~\ref{lem02} and~\ref{lem02a}, for proving that the conjecture holds for 
$C_3$-free planar graphs it suffices to prove only $C_3$-free planar graphs 
 of order $\D t$ where $t\geq 3$ is a positive integer.
\end{obs}

\begin{lem}\label{lem03}~\cite{ZY98}
	Let $G$ be a graph of order $mt$ with chromatic number $\chi $ such that $\chi \leq m,$
	where $t$ is an integer.
		If $e(G) \leq (m-1)t,$ then $G$ is equitably $m$-colorable.
\end{lem}

\begin{lem}\label{lem04} Suppose $G$ is a $C_3$-free planar graph
with $\D(G) = \D.$ If $G$ has an independent $s$-set $V'$ and there
exists $U \subseteq V(G) - V'$ such that $|U| > s(1+\D )/2$
and $e(u, V') \geq 1$ for all $u \in U,$ then $U$ contains two
nonadjacent vertices $\alpha$ and $\beta$ which are adjacent to
exactly one and the same vertex $\gamma \in V'.$
\end{lem}

\noindent \textbf{Proof.}
Let $ U_1$ consist of vertices in $U$ with exactly one neighbor in $V'.$
If $r=|U_1|,$ then $r+2(|U| - r)\leq \D s$ which implies $r \geq 2|U|- \D s >s.$
Consequently, $V'$ contains a vertex $\gamma$ which has at least two neighbors in $U_1.$
Since $G$ is $C_3$-free, this two neighbors are not adjacent. 
Thus $U_1$ contains  two nonadjacent vertices $\alpha$ and $\beta$
which are adjacent to exactly one and the same vertex $\gamma \in V'.$ \qed

\begin{lem} \label{lem05} ~\cite{Nak12}
If a graph $G$ has an independent $s$-set $V'$ and there exists
$U \subseteq V(G) - V'$  such that  $e(u,V') \geq 1 $ for all $u \in U,$
and $e(G[U])+e(V',U) < 2 |U|  -  s,$ then
 $U$ contains two nonadjacent vertices $\alpha$ and $\beta$ which are
adjacent to exactly one and the same vertex $\gamma \in V'.$
\end{lem}

\noindent \textbf{Notation.} {\em Let $q_{m,\D,t}$ denote the maximum
number not exceeding $2mt - 4$ such that each $C_3$-free planar graph of order
$mt,$ where $t$ is an integer, is equitably $m$-colorable if
it has maximum degree at most $\D$ and size at most $q_{m,\D,t}.$}

The next Lemma is similar to that in~\cite{Nak12} except that we use $V_1$
instead of  $V'_1$ which is erratum. Nevertheless later arguments in~\cite{Nak12} stand correct.

\begin{lem}\label{lem07a}
Let $G$ be an edge-minimal $C_3$-free planar graph that is not equitably $m$-colorable
with order $mt,$ where $t$ is an integer, and maximum
degree at most $\D.$ If $e(G) \leq
(r+1)(m-r)t-t+2+q_{r,\D,t},$ then $B$ contains two nonadjacent
vertices $\alpha$ and $\beta$ which are adjacent to exactly one and
the same vertex $\gamma \in V_1.$
 \end{lem}

\noindent \textbf{Proof.}
If $e(G[A']) \leq q_{r,\D,t},$ then $G[A']$ is equitably $r$-colorable.
Consequently, $G$ is equitably $m$-colorable.
So we suppose $e(G[A']) \geq q_{r,\D,t}+1.$
By Observation~\ref{obs01}, $e(A' - V'_1, B') \geq (r-1)(m-r)t.$
Note that $e(G[B])=e(G[B']), e(V_1,B)= e(V'_1,B')+1.$
So $e(G[B]) + e(V_1,B) =e(G[B']) + e(V'_1,B')+1
= e(G) - e(G[A']) -e(A' - V'_1, B')+1 < 2mt-2rt -t +3 = 2|B| - |V_1|.$
By Lemma~\ref{lem05}, $B$ contains two nonadjacent vertices $\alpha$ and $\beta$
which are adjacent to exactly one and the same vertex $\gamma \in V_1.$ \qed

\begin{lem}\label{lem07b}
Let $G$ be an edge-minimal $C_3$-free planar graph that is not equitably $m$-colorable
with order $mt,$ where $t$ is an integer, and maximum
degree at most $\D.$ If $B$ contains two
nonadjacent vertices $\alpha$ and $\beta$ which are adjacent to
exactly one and the same vertex $\gamma \in V_1,$ then $e(G) \geq
r(m-r)t+q_{r,\D,t}+q_{m-r,\D,t}-\D+4.$
\end{lem}
\noindent \textbf{Proof.} Suppose $e(G) \leq r(m-r)t+q_{r,\D,t}+q_{m-r,\D,t}-\Delta+3.$
If $e(G[A']) \leq q_{r,\D,t},$ then $G[A']$ is equitably $r$-colorable.
Consequently, $G$ is equitably $m$-colorable.
So we suppose $e(G[A']) \geq q_{r,\D,t}+1.$
This with Observation~\ref{obs01} implies $e(G[A'])+ e(A,B') \geq q_{r,\D,t}+1+r(m-r)t.$
Note that $e(G[A'])+ e(A,B') = e(G[A])+e(A,B).$
Let $A_1 =  (A - \{ \gamma \}) \cup \{ \alpha, \beta \}$ and
$B_1 = (B \cup \{ \gamma \}) - \{ \alpha, \beta \}.$
Then $e(G[A_1])+e(A_1, B_1)  \geq  e(G[A])+e(A,B) - \D +2 \geq q_{r,\D,t}+1+r(m-r)t - \D +2.$
So $e(G[B_1]) = e(G) - e(G[A_1])+e(A_1, B_1) \leq q_{m-r,\D,t}$ which implies
$G[B_1]$ is equitably $(m-r)$-colorable.
Combining with $(V_1 - \{\gamma\}) \cup \{\alpha, \beta\},V_2,\ldots,V_r,$ we have $G$
equitably $m$-colorable which is a contradiction. \qed

\begin{cor}\label{cor01}
Let $G$ be an edge-minimal $C_3$-free planar graph that is not equitably $m$-colorable
with order $mt,$ where $t$ is an integer, and maximum
degree at most $\D.$ Then $e(G) \geq
r(m-r)t+q_{r,\D,t}+q_{m-r,\D,t}-\D+4$ if one of the following conditions
are satisfied: \\ (i) $(m-r)t +1> (t-1)(1+\D )/2;$\\
(ii) $e(G) \leq (r+1)(m-r)t-t+2+q_{r,\D,t}.$
\end{cor}

\noindent \textbf{Proof.} This is a direct consequence of Lemmas
\ref{lem04},~\ref{lem07a}, and~\ref{lem07b}. \qed

Now we are ready to work on $C_3$-free planar graphs. 

\begin{lem}\label{lem08}
(i) $q_{1,\D,t}=0.$ \quad (ii) $q_{2,\D,t}\geq 3$ for $t\geq 3.$
\quad (iii) $q_{3,\D,t} \geq 2t.$
\end{lem}
\noindent \textbf{Proof.} (i) and (ii) are obvious. (iii) is the result of Theorem~\ref{Grot} 
and Lemma~\ref{lem03}. \qed

\begin{lem}\label{lem16}
$q_{4,\Delta,t} \geq \min\{q_{3,\Delta,t}+3t+3-\D, 4t-\D+9 \}$ 
for $\D \geq 5$ and $t\geq 3.$
\end{lem}

\noindent \textbf{Proof.} Condider $\D \geq 5$ and $t \geq  3.$
Suppose $G'$ is a $C_3$-free planar graph with maximum degree at most $\D$ and
$e(G') \leq  \min\{q_{3,\Delta,t}+3t+3-\D, 4t-\D+9 \}$
 but $G'$ is not equitably $4$-colorable.
Let $G \subseteq G'$ be an edge-minimal graph that is not equitably $4$-colorable.
From Table~\ref{tab01}, $e(G) > e(G').$ This contradiction completes the proof.

\begin{table}[h]
\centering
    \begin{tabular}{|c|c|c|}\hline
         $r$   & lower bounds on size  & Reasons  \\ \hline
         3  & $q_{3,\Delta,t}+3t+3-\D$ or $q_{3,\Delta,t}+3t+2$ & Corollary~\ref{cor01}(ii), Lemma~\ref{lem08}\\
         2  & $4t-\D+9$ or $5t+5$ & Corollary~\ref{cor01}(ii), Lemma~\ref{lem08}\\
         1  & $q_{3,\Delta,t}+3t+3-\D$ or $5t+2$ &   Corollary~\ref{cor01}(ii), Lemma~\ref{lem08}\\   \hline
    \end{tabular}
\caption{Lower bounds on size of $G$ in the proof of Lemma~\ref{lem16}}\label{tab01}
\end{table}
\qed

\begin{lem}\label{lem17}
$q_{5,\Delta,t} \geq \min\{q_{3,\Delta,t}+6t+6-\D, q_{4,\Delta,t}+4t+3-\D , 7t+2 \}$ 
for $\D \geq 5$ and $t\geq 3.$
\end{lem}

\noindent \textbf{Proof.} Use Table~\ref{tab02} for an argument similar to the proof of Lemma~\ref{lem16}. 

\begin{table}[h]
\centering
    \begin{tabular}{|c|c|c|}\hline
         $r$   & lower bounds on size  & Reasons  \\ \hline
         3  & $q_{3,\Delta,t}+6t+6-\D$ or $q_{3,\Delta,t}+7t+2$ & Corollary~\ref{cor01}(ii), Lemma~\ref{lem08}\\
         2  & $q_{3,\Delta,t}+6t+6-\D$  or $8t+5$ & Corollary~\ref{cor01}(ii), Lemma~\ref{lem08}\\
         1  & $q_{4,\Delta,t}+4t+3-\D$  & Corollary~\ref{cor01}(i), Lemma~\ref{lem08}\\   \hline
    \end{tabular}
\caption{Lower bounds on size of $G$ in the proof of Lemma~\ref{lem17}}\label{tab02}
\end{table}
\qed

\begin{cor}\label{cor05} (1) $q_{4,6,t}$ is at least $5t-3$ and $4t+3$ for $t$ at least $3$ and $6,$ respectively.\\
\noindent (2) $q_{4,7,t}$ is at least $5t-4$ and $4t+2$ for $t$ at least $3$ and $6,$ respectively.\\
\noindent (3) $q_{5,6,t}$ is at least $9t-6$ and $8t$ for $t$ at least $3$ and $6,$ respectively.\\
\noindent (4) $q_{5,7,t}$ is at least $9t-8$ and $8t-2$ for $t$ at least $3$ and $6,$ respectively.
\end{cor}

\noindent \textbf{Proof.} The results can be calculated directly from Lemmas~\ref{lem08} to~\ref{lem17}. \qed

\begin{cor}\label{cort1}
Each $C_3$-free planar graph $G$ with maximum degree at most $7$ 
and $|G|\geq 18$ has an equitable $6$-coloring. 
Moreover, each $C_3$-free planar graph $G$ with maximum degree $6$  
has an equitable $6$-coloring. 
\end{cor}
\noindent \textbf{Proof.} 
Let $G$ be an edge-minimal $C_3$-free planar graph 
that is not equitably $\D$-colorable with $|G|=6t,$ where $t$ is an integer at least $3,$  
and maximum degree at most $7.$\\
\indent Consider the case $r = 3.$ 
By Corollaries~\ref{cor01}(ii) and~\ref{cor05}, 
$e(G) > \min \{2q_{3,\Delta,t}+9t+3-\D,q_{3,\Delta,t}+11t+2\} \geq 13t-4 \geq 12t-4.$\\
\indent Consider the case $r = 2.$ 
By Corollary~\ref{cor01}(i), $e(G)>q_{4,\Delta,t}+8t+6-\D.$ 
It follows from Corollary~\ref{cor05} that $e(G)>\min\{13t-5, 12 t+1\}\geq 12 t-4$ for $t\geq 3.$ \\
\indent Consider the case $r = 1.$ 
We have $e(B', V_1) \geq 5t$ by Observation~\ref{obs01}. 
But $y$ has at most $\D -1$ neighbors in $B'$ because $xy \in E(G),$ 
so $(t-1)\D - 1 \geq e(B',V_1).$ 
Consequently, $(t-1)\D - 1 \geq 5t.$ That is $t \geq 4$ when $\D \leq 7.$  
By Corollary~\ref{cor01} (i), $e(G) > q_{5,\D,t}+5t-4.$ 
Using Corollary~\ref{cor05}, we have $e(G)>\min\{14t-12, 13t-6\}.$ 
 It follows from $t\geq 4$ that $e(G)> 12t-4.$   \\
\indent Since we have contradiction for all cases, the counterexample is impossible. 
Use Lemma~\ref{lem02} to complete the first part of the proof. \\
\indent Observation~\ref{obs02} implies each $C_3$-free planar graph $G$ with 
maximum degree $6$ has an equitable $6$-coloring.  \qed

Note that a graph $G$ in Corollary~\ref{cort1} 
has an equitable $m$-coloring with $m< \D(G).$ 

\begin{lem}\label{lem18}
Each $C_3$-free planar graph $G$ with maximum degree at most $7$ 
 has an equitable $7$-coloring. 
\end{lem}

\noindent \textbf{Proof.} Use Table~\ref{tab03} for an argument similar to the proof of Lemma~\ref{lem16}. 

\begin{table}[h]
\centering
    \begin{tabular}{|c|c|c|}\hline
         $r$   & lower bounds on size  & Reasons  \\ \hline
         3  & $q_{3,\Delta,t}+12t+q_{4,\Delta,t}+3-\D$ &Corollary~\ref{cor01}(i), Lemma~\ref{lem08}\\
         2  & $q_{5,\Delta,t}+10t+6-\D$  & Corollary~\ref{cor01}(i), Lemma~\ref{lem08}\\
         1  & $q_{6,\Delta,t}+6t+3-\D$  & Corollary~\ref{cor01}(i), Lemma~\ref{lem08}\\   \hline
    \end{tabular}
\caption{Lower bounds on size of $G$ in the proof of Lemma~\ref{lem18}}\label{tab03}
\end{table}
\qed

\indent Using Corollary~\ref{cor05}, and $q_{6,\Delta,t}=12t-4$ from Corollary~\ref{cort1},  
we have $e(G)> 14 t -4$ for each case of $r,$ which is a contradiction. 
Thus the counterexample is impossible. 
Use Observation~\ref{obs02} to complete the proof. \qed 

\begin{Th}\label{cort2}
Each $C_3$-free planar graph $G$ with maximum degree $\D \geq 6$ 
 has an equitable $\D$-coloring. 
\end{Th}
\noindent \textbf{Proof.} 
Zhu and Bu~\cite{ZhuBu08} proved that the theorem holds for every $C_3$- free planar graph
with maximum degree at least $8.$  
Use Corollary~\ref{cort1} and Lemma~\ref{lem18} to complete the proof. \qed

Next, we show that the conjecture holds also for a planar graph of maximum degree $5$ 
if we restrict the girth to be at least $6.$  

\begin{cor}\label{cort3}
Each planar graph $G$ of girth at least $6$ with maximum degree at most $6$ 
and $|G|\geq 15$ has an equitable $5$-coloring. 
Moreover, each planar graph $G$ with girth at least $6$ and 
maximum degree $\D \geq 5$ has an equitable $\D$-coloring. 
\end{cor}
\noindent \textbf{Proof.} 
Let $G$ be an edge-minimal planar graph of girth at least $6$ 
that is not equitably $\D$-colorable with $|G|=5t,$ where $t$ is an integer at least $3,$  
and maximum degree at most $6.$\\
 \indent Then for $t \geq 3,$ we have $e(G) \leq (15/2)t -3$ from Lemma~\ref{edge}, 
and $e(G) > \min\{9t-6, 8t\}$ from Corollary~\ref{cor05} which leads to a contradiction. 
Thus the counterexample is impossilble. 
Use Lemma~\ref{lem02} to complete the first part of the proof. \\
\indent Observation~\ref{obs02} implies each planar graph $G$ with girth at least $6$ and 
maximum degree $5$ has an equitable $5$-coloring. 
Use Theorem~\ref{cort2} to complete the proof. \qed

\section{Results on Planar Graphs without $C_4$} 

First we introduce the result by Tan~\cite{Tan10}. 

\begin{lem}\label{edge2} 
If a planar graph $G$ of order $n$ does not contain $C_4,$ then $e(G) \leq (15/7) n - (30/7)$ and 
$\delta(G) \leq 4.$ 
\end{lem} 
The proof of Lemma~\ref{edge2} by Tan is presented here for convenience of readers. \\

\noindent \textbf{Proof.} 
Let $f$ and $f_i$ denote the number of faces 
and the number of faces of length $i,$ respectively. 
We need only to consider the case that $G$ is connected. 
A graph $G$ cannot contain two $C_3$ that share the same edge 
since $G$ does not contain $C_4.$ It follows that $3f_3 \leq e(G).$ \\
\indent Consider  $5f-2f_3 = 5(f_3+f_5 +\cdots+f_n) - 2f_3 \leq 3f_3  +5f_5 +\cdots+nf_n 
=\sum_{1 \leq i \leq n}if_i = 2e(G).$ 
Thus $ f \leq (8/15)e(G).$ 
Using Euler's formula, we have  $e(G) \leq (15/7) n - (30/7).$
The result about minimum degree follows from Handshaking Lemma. 
 \qed

From Lemma~\ref{edge2}, each edge-minimal counterexample graph has $1 \leq r \leq 4.$ 
The following tools in this section are quite similar to that of the previous section. 
Thus we omit the proofs of them. 

\begin{lem}\label{lem402}
Let $m$ be a fixed integer with $m \geq 1.$  Suppose that any 
planar graph without $C_4$ of order $mt$ with maximum degree at most $\D$ is
equitably $m$-colorable for any integer $t \geq k.$ Then any
planar graph without $C_4$ of order at least $kt$ and maximum degree
at most $\D$ is also equitably $m$-colorable.
\end{lem}

\begin{obs}\label{obs402}  By Lemmas~\ref{lem02a} and~\ref{lem402}, 
for proving that the conjecture holds for 
planar graphs without $C_4$ it suffices to prove only planar graphs without $C_4$
 of order $\D t$ where $t\geq 3$ is a positive integer.
\end{obs}

\begin{lem}\label{lem404} Suppose $G$ is a planar graph without $C_4$ 
with $\D(G) = \D.$ If $G$ has an independent $s$-set $V'$ and there
exists $U \subseteq V(G) - V'$ such that $|U| > s(2+\D )/2$
and $e(u, V') \geq 1$ for all $u \in U,$ then $U$ contains two
nonadjacent vertices $\alpha$ and $\beta$ which are adjacent to
exactly one and the same vertex $\gamma \in V'.$
\end{lem}

\noindent \textbf{Notation.} {\em Let $p_{m,\D,t}$ denote the maximum
number not exceeding $(15/7)mt - (30/7)$ such that each planar graph without $C_4$ 
of order $mt,$ where $t$ is an integer, is equitably $m$-colorable if
it has maximum degree at most $\D$ and size at most $p_{m,\D,t}.$}

\begin{lem}\label{lem407a}
Let $G$ be an edge-minimal planar graph without $C_4$ 
that is not equitably $m$-colorable 
with order $mt,$ where $t$ is an integer, and maximum
degree at most $\D.$ If $e(G) \leq
(r+1)(m-r)t-t+2+p_{r,\D,t},$ then $B$ contains two nonadjacent
vertices $\alpha$ and $\beta$ which are adjacent to exactly one and
the same vertex $\gamma \in V_1.$
 \end{lem}

\begin{lem}\label{lem407b}
Let $G$ be an edge-minimal planar graph without $C_4$ that is not equitably $m$-colorable
with order $mt,$ where $t$ is an integer, and maximum
degree at most $\D.$ If $B$ contains two
nonadjacent vertices $\alpha$ and $\beta$ which are adjacent to
exactly one and the same vertex $\gamma \in V_1,$ then $e(G) \geq
r(m-r)t+p_{r,\D,t}+p_{m-r,\D,t}-\D+4.$
\end{lem}

\begin{cor}\label{cor401}
Let $G$ be an edge-minimal planar graph without $C_4$ 
that is not equitably $m$-colorable 
with order $mt,$ where $t$ is an integer, and maximum 
degree at most $\D.$ Then $e(G) \geq
r(m-r)t+p_{r,\D,t}+p_{m-r,\D,t}-\D+4$ if one of the following conditions
are satisfied:\\ (i) $(m-r)t +1> (t-1)(2+\D )/2;$\\
(ii) $e(G) \leq (r+1)(m-r)t-t+2+p_{r,\D,t}.$
\end{cor}

Now we are ready to work on planar graphs without $C_4.$ 

\begin{lem}\label{lem408}
(i) $p_{1,\D,t}=0.$ \quad (ii) $p_{2,\D,t} = 2.$ 
\quad (iii) $p_{3,\D,t} \geq 6$ for $t \geq  3.$ \quad (iv) $p_{4,\D,t} \geq 3t.$
\end{lem}
\noindent \textbf{Proof.} (i), (ii) and (iii) are obvious. (iv) is the result of  Lemma~\ref{lem03}. \qed

\begin{lem}\label{lemp5}
$p_{5,\Delta,t} \geq \min\{p_{4,\Delta,t}+16t+3-\D, 6t+11-\D, 7t+2 \}$ 
for $\D \geq 8$ and $t\geq 3.$
\end{lem}

\noindent \textbf{Proof.} Use Table~\ref{tabp5} for an argument similar to the proof of Lemma~\ref{lem16}. 

\begin{table}[h]
\centering
    \begin{tabular}{|c|c|c|}\hline
         $r$   & lower bounds on size  & Reasons  \\ \hline
         4  & $p_{4,\Delta,t}+16t+3-\D$ or $p_{4,\Delta,t}+4t+2$  & Corollary~\ref{cor401}(ii), Lemma~\ref{lem408}\\
         3  & $6t+11-\D$ or $7t+8$ & Corollary~\ref{cor401}(ii), Lemma~\ref{lem408}\\
         2  & $6t+11-\D$ or $8t+4$ & Corollary~\ref{cor401}(ii), Lemma~\ref{lem408}\\
         1  & $p_{4,\Delta,t}+4t+3-\D$  or $7t+2$ & Corollary~\ref{cor401}(ii), Lemma~\ref{lem408}\\   \hline
    \end{tabular}
\caption{Lower bounds on size of $G$ in the proof of Lemma~\ref{lemp5}}\label{tabp5}
\end{table}
\qed

\begin{lem}\label{lemp6}
$p_{6,\Delta,t} \geq \min\{p_{4,\Delta,t}+8t+5-\D,9t+15-\D,11t+4, p_{5,\Delta,t}+5t+3-\D \}$ 
for $\D \geq 8$ and $t\geq 3.$
\end{lem}

\noindent \textbf{Proof.} Use Table~\ref{tabp6} for an argument similar to the proof of Lemma~\ref{lem16}. 

\begin{table}[h]
\centering
    \begin{tabular}{|c|c|c|}\hline
         $r$   & lower bounds on size  & Reasons  \\ \hline
         4  & $p_{4,\Delta,t}+8t+5-\D$ or $p_{4,\Delta,t}+9t+2$  & Corollary~\ref{cor401}(ii), Lemma~\ref{lem408}\\
         3  & $9t+15-\D$ or $11t+8$ & Corollary~\ref{cor401}(ii), Lemma~\ref{lem408}\\
         2  & $p_{4,\Delta,t}+8t+5-\D$ or $11t+4$ & Corollary~\ref{cor401}(ii), Lemma~\ref{lem408}\\
         1  & $p_{5,\Delta,t}+5t+3-\D$  & Corollary~\ref{cor401}(i), Lemma~\ref{lem408}\\   \hline
    \end{tabular}
\caption{Lower bounds on size of $G$ in the proof of Lemma~\ref{lemp6}}\label{tabp6}
\end{table}
\qed

\begin{cor}\label{cor405} (1) $p_{5,8,t}$ is at least $7t-5$ and $6t+3$ for $t$ at least $3$ and $8,$ respectively.\\
\noindent (2) $p_{6,8,t}$ is at least $12t-10$ and $9t+7$ for $t$ at least $3$ and $6,$ respectively.\\
\noindent (3) $p_{7,8,t}$ is at least $18t-15$ and $15t+1$ for $t$ at least $3$ and $6,$ respectively.
\end{cor}

\noindent \textbf{Proof.} The results can be calculated directly from Lemmas~\ref{lem408} to~\ref{lemp6}. \qed

\begin{cor}\label{cortp1}
Each planar graph $G$ without $C_4$ with maximum degree at most $8$ 
and $|G|\geq 21$ has an equitable $7$-coloring. 
Moreover, each planar graph $G$ without $C_4$ with maximum degree  $7$  
has an equitable $7$-coloring.  
\end{cor}

\noindent \textbf{Proof.} 
Let $G$ be an edge-minimal planar graph without $C_4$  
that is not equitably $\D$-colorable with $|G|=7t,$ where $t$ is an integer at least $3,$  
and maximum degree at most $8.$\\
\indent Consider the case $r = 4.$ 
By Corollaries~\ref{cor401}(ii) and~\ref{cor405}, 
$e(G) > \min \{p_{4,\Delta,t}+12t+p_{3,\Delta,t}+3-\D,p_{4,\Delta,t}+14t+2\} \geq 15t+1 \geq 15t - (30/7)$ 
for $t\geq 3.$ \\
\indent Consider the case $r = 3.$ 
By Corollaries~\ref{cor401}(ii) and~\ref{cor405}, 
$e(G) > \min \{p_{3,\Delta,t}+12t+p_{4,\Delta,t}+3-\D,p_{3,\Delta,t}+15t+2\} \geq 15t+1 \geq 15t - (30/7)$ 
for $t\geq 3.$ \\
\indent Consider the case $r = 2.$ 
By Corollaries~\ref{cor401}(i) and~\ref{cor405}, 
$e(G) >   10t+p_{5,\Delta,t}+3-\D  \geq 15t - (30/7)$ for $t \geq 3.$ \\
\indent Consider the case $r = 1.$ 
We have $e(B', V_1) \geq 6t$ by Observation~\ref{obs01}. 
But $y$ has at most $\D -1$ neighbors in $B'$ because $xy \in E(G),$ 
so $(t-1)\D - 1 \geq e(B',V_1).$ 
Consequently, $(t-1)\D - 1 \geq 5t.$ That is $t \geq 4.5$ when $\D \leq 8.$  
By Corollary~\ref{cor401} (i), $e(G) > p_{6,\D,t}+6t-5.$ 
Using Corollary~\ref{cor405}, we have $e(G)>\min\{18t-15, 15t+1\}.$  
 It follows from $t\geq 4.5$ that $e(G)> 15t - (30/7).$   \\
\indent Since we have contradiction for all cases, the counterexample is impossible. 
Use Lemma~\ref{lem402} to complete the first part of the proof. \\
\indent Observation~\ref{obs402} implies each planar graph $G$ 
without $C_4$ with maximum degree $7$ has an equitable $7$-coloring.  \qed

\begin{lem}\label{lem418}
Each planar graph $G$ without $C_4$ with maximum degree at most $8$ 
 has an equitable $8$-coloring. 
\end{lem}

\noindent \textbf{Proof.} Use Table~\ref{tabp03} for an argument similar to the proof of Lemma~\ref{lem16}. 

\begin{table}[h]
\centering
    \begin{tabular}{|c|c|c|}\hline
         $r$   & lower bounds on size  & Reasons  \\ \hline
         4  & $p_{4,\Delta,t}+16t+p_{4,\Delta,t}+3-\D$ or $p_{4,\Delta,t}+19t+2$  & Corollary~\ref{cor401}(ii), Lemma~\ref{lem408}\\
         3  & $15t+p_{5,\Delta,t}+9-\D$  & Corollary~\ref{cor401}(i), Lemma~\ref{lem408}\\
         2  & $p_{6,\Delta,t}+12t+5-\D$ & Corollary~\ref{cor401}(i), Lemma~\ref{lem408}\\
         1  & $p_{7,\Delta,t}+7t+3-\D$  & Corollary~\ref{cor401}(i), Lemma~\ref{lem408}\\   \hline
    \end{tabular}
\caption{Lower bounds on size of $G$ in the proof of Lemma~\ref{lem418}}\label{tabp03}
\end{table}
\qed

\indent Using Corollary~\ref{cor405}, and $q_{7,\Delta,t}=15t-4$ from Corollary~\ref{cortp1},  
we have $e(G)> (120/7) t - (30/7)$ for each case of $r,$ which is a contradiction. 
Thus the counterexample is impossible. 
Use Observation~\ref{obs402} to complete the proof. \qed 

\begin{Th}\label{cortp2}
Each  planar graph $G$ without $C_4$ with maximum degree $\D \geq 7$ 
 has an equitable $\D$-coloring. 
\end{Th}
\noindent \textbf{Proof.} 
Nakprasit~\cite{Nak12} proved that the theorem holds for every planar graph
with maximum degree at least $9.$  
Use Corollary~\ref{cortp1} and Lemma~\ref{lem418} to complete the proof. \qed

\section{Acknowledgement}
The first author was supported by Research Promotion Fund, Khon Kaen University, Fiscal year 2011.

\end {document}